\documentclass{article}
\usepackage{amssymb, amsthm}
\usepackage{pifont,mathrsfs}
\usepackage{array,lastpage}
\usepackage{enumerate,xspace,pifont,shadow}
\usepackage{mathrsfs}
\usepackage{graphicx}
\usepackage{multimedia}
\usepackage{calrsfs}
\usepackage[all]{xy}

\DeclareSymbolFont{AMSb}{U}{msb}{m}{n}
\DeclareMathSymbol{\Z}{\mathbin}{AMSb}{"5A}
\DeclareMathSymbol{\R}{\mathbin}{AMSb}{"52}
\DeclareMathSymbol{\N}{\mathbin}{AMSb}{"4E}
\DeclareMathSymbol{\Q}{\mathbin}{AMSb}{"51}

\newcommand{\tp}{\textup{tp}}
\newcommand{\stp}{\textup{stp}}
\newcommand{\dcl}{\textup{dcl}}
\newcommand{\acl}{\textup{acl}}
\newcommand{\Th}{\textup{Th}}
\newcommand{\dom}{\textup{dom}}
\newcommand{\F}{\textup{\textbf{F}}}
\newcommand{\cf}{\textup{cf}}

\newcommand{\Inf}{\textup{Inf}}

\newcommand{\Aut}{\textup{Aut}}
\newcommand{\ran}{\textup{ran}}
\newcommand{\Av}{\textup{Av}}

\def\Ind{\setbox0=\hbox{$x$}\kern\wd0\hbox to 0pt{\hss$\mid$\hss}
\lower.9\ht0\hbox to 0pt{\hss$\smile$\hss}\kern\wd0}

\def\Notind{\setbox0=\hbox{$x$}\kern\wd0\hbox to 0pt{\mathchardef
\nn=12854\hss$\nn$\kern1.4\wd0\hss}\hbox to
0pt{\hss$\mid$\hss}\lower.9\ht0 \hbox to
0pt{\hss$\smile$\hss}\kern\wd0}

\def\ind{\mathop{\mathpalette\Ind{}}}

\def\nind{\mathop{\mathpalette\Notind{}}}

\newtheorem{thm}{Theorem}[section]
\newtheorem{lem}[thm]{Lemma}
\newtheorem{cor}[thm]{Corollary}
\newtheorem{prop}[thm]{Proposition}
\newtheorem{conj}[thm]{Conjecture}
\newtheorem{fact}[thm]{Fact}

\newtheorem{quest}[thm]{Question}

\theoremstyle{definition}
\newtheorem{definition}[thm]{Definition}

\theoremstyle{remark}
\newtheorem{remark}[thm]{Remark}

\theoremstyle{remark}
\newtheorem{example}[thm]{Example}

\theoremstyle{remark}
\newtheorem{claim}[thm]{Claim}

\theoremstyle{remark}
\newtheorem{notation}[thm]{Notation}

\begin{document}

\bibliographystyle{plain}

\title{When does elementary bi-embeddability imply isomorphism?}
\author{John Goodrick}

\maketitle

\begin{abstract}
A first-order theory has the \emph{Schr\"{o}der-Bernstein property} if any two of its models that are elementarily bi-embeddable are isomorphic.  We prove that if a countable theory $T$ has the Schr\"{o}der-Bernstein property then it is classifiable (it is superstable and has NDOP and NOTOP) and satisfies a slightly stronger condition than nonmultidimensionality, namely: there cannot be $M \models T$, $p \in S(M)$, and $f \in \Aut(M)$ such that for every $i < j < \omega$, $f^i(p) \bot f^j(p)$.  We also make some conjectures about how the class of theories with the Schr\"{o}der-Bernstein property can be characterized.

\end{abstract}

\section{History and summary of results}

This paper is about the following property of a first-order theory $T$:

\begin{definition}
A theory $T$ has the \emph{Schr\"{o}der-Bernstein property}, or the \emph{SB~property} for short, if it is the case that whenever $M, N \models T$, $M$ is elementarily embeddable into $N$, and $N$ is elementarily embeddable into $M$, then $M$ and $N$ are isomorphic.
\end{definition}

Our main goal is to find some nice model-theoretic characterization of the class of complete theories with the SB~property.  This property was first studied in the 1980's by Nurmagambetov in \cite{nur2} and \cite{nur1}, mainly within the class of $\omega$-stable theories.  In \cite{nur1}, he showed:

\begin{thm}
\label{ttrans}
If $T$ is $\omega$-stable, then $T$ has the SB~property if and only if $T$ is nonmultidimensional.
\end{thm}

The SB~property for more general classes of theories is also related to concepts from Shelah's classification theory.  For instance, we will prove the following below:

\begin{thm}
(Corollary~\ref{class_SB}) If $T$ is a countable theory and $T$ has the SB~property then $T$ is classifiable.
\end{thm}

In the theorem above, ``classifiable'' is to be understood in a Shelahian sense: it means roughly that every model can be characterized by a tree of cardinal invariants.  For our purposes, we could also define $T$ to be classifiable if it falls on the lower side of the Main Gap expounded in \cite{bible}.  In the case where $T$ is countable, $T$ is classifiable if and only if it is superstable and also satisfies the more technical conditions of having NDOP and NOTOP, and a further hypothesis of shallowness is needed for $T$ to have ``few models'' (see \cite{bible} for details).

But not all classifiable theories have the SB~property:

\begin{thm}
(Corollary~\ref{superstableSB}) A superstable theory $T$ does not have the SB~property if $T$ has a \emph{nomadic} type, that is, a type $p \in S(M)$ such that there is an automorphism $f \in \Aut(M)$ for which the types $\left\{f^n(p) : n \in \N \right\}$ are pairwise orthogonal.
\end{thm}

In the next section we discuss some motivating conjectures about the SB~property.  In section~3 we give examples of theories and check the SB~property for each of them; we think that these examples give a feeling for the kinds of difficulties that can arise in showing that more general classes of theories do, or do not, have the SB~property.  In section~4 we prove a strong negation of the SB~property for any theory in which arbitrarily long orders are definable in $L_{\infty, \omega}$ (for instance, unstable theories, or stable theories with OTOP).  In section~5, we discuss what we know about the stable case: we use prime model machinery and forking calculus to prove that the SB~property fails for non-superstable $T$ and stable $T$ with regular nomadic types, and we prove that the SB~property holds for perfectly~trivial superstable theories with no nomadic types.

Most of our notation is standard for model~theory and consistent with that in \cite{bible}.  For a general background in the stability theory used in this paper, the reader is referred to \cite{baldwin}, \cite{pillay}, and \cite{poizat}.

Throughout this paper, $T$ denotes a complete first-order theory.  Note that $T$ is not assumed to be countable unless we explicitly say so.

We will assume for convenience that all the sets in this paper live in a ``monster model'' $\mathfrak{C}$ that is $\kappa$-saturated and strongly $\kappa$-homogeneous for some uncountable cardinal $\kappa$ that is larger than the size of any of the other sets we will consider.  By standard arguments, such a model exists for any theory.

\section{Some conjectures}

We list here our main conjectures about the Schr\"{o}der-Bernstein property.  All of the results in this paper came from attempts to prove or refute conjectures from this section.

Our main goal is to prove a characterization of the Schr\"{o}der-Bernstein property along the following lines:

\begin{conj}
\label{main_SB_conj}
The complete theory $T$ has the SB~property if and only if it has both of the following properties:

1. $T$ is classifiable,

2. $T$ has no \emph{wandering types}: that is, there is no $M \models T$, $p \in S(M)$, and $f \in \Aut(M)$ such that for every $i < j < \omega$, $f^i(p) \bot^a f^j(p)$.

\end{conj}

Here is roughly how we expect that Conjecture~\ref{main_SB_conj} could be proved: first, we know that a theory with the SB~property must be classifiable, at least in the countable case (and unclassifiable theories in uncountable languages still ought to be ``wild'' enough to contain counterexamples to SB).  We also know that such a theory cannot contain any nomadic types (``nomadic'' has the same definition as ``wandering'' except that we replace $\bot^a$ by $\bot$).  If a classifiable theory contains a type that is wandering but not nomadic, then we expect that it is nonorthogonal to a definable regular locally~modular group -- maybe this is provable by arguments similar to those in \cite{almost_orthogonality}.  Hopefully this will allow us to reduce to proving that a classifiable, locally~modular abelian group has the SB~property if and only if it has no wandering types (though currently there seem to be difficult obstacles to carrying out this reduction).  Some special cases of this conjecture for locally~modular groups have been checked already in \cite{ab_groups} and \cite{my_thesis}.  For instance, in \cite{ab_groups} it is shown that for any abelian group $G$, $\Th(G; +)$ has the SB~property if and only if $G$ is superstable and has no wandering types; and if $T$ is a weakly minimal theory of an abelian group, possibly in a language larger than $\left\{+\right\}$, it is shown without too much difficulty in \cite{my_thesis} that if $T$ has no wandering types then $T$ has the SB~property.

We also conjecture that if the SB~property fails on some definable set, then it fails for the entire theory:

\begin{conj}
\label{lifting-SB-conj}
Suppose $T$ is a complete theory $M \models T$, $\varphi(\bar{x})$ is a formula in $T$, and $T'$ is the theory of the structure $\varphi(M)$ with all the induced definable structure from $M$.  If $T'$ does not have the SB~property, then neither does $T$.
\end{conj}

Given a formula $\varphi(\bar{x})$ in $T$ whose theory $T'$ (as above) does not have the SB~property, the na\"{i}ve way to approach Conjecture~\ref{lifting-SB-conj} is to take two models $M'$ and $N'$ of $T'$ that are bi-embeddable but nonisomorphic and ``lift'' them to bi-embeddable models $M$ and $N$ of $T$ such that $\varphi(M) = M'$ and $\varphi(N) = N'$, and then $M \ncong N$ since they have nonisomorphic interpretations of $\varphi(\bar{x})$.  However, it is nontrivial to get a even a single model $M \models T$ such that $\varphi(M) = M'$; $T$ is stable and \emph{countable}, this is possible using ``$\ell$-constructions'' (see Remark~4.3 of \cite{rem_nmd_sstable}).  Furthermore, it is unclear how (or whether) one can pick $\ell$-constructions over $M'$ and $N'$ that result in bi-embeddable models of $T$, since the models $M$ and $N$ are not ``prime'' in the usual sense over $M'$ and $N'$.

Another question is how sensitive the SB~property is to the addition of named constants to the language of $T$.  In fact naming a single constant can cause a theory without the SB~property to gain the SB~property -- see Example~\ref{single_const}.  But we do not think that naming constants can cause a theory to \emph{lose} the SB~property:

\begin{conj}
1. There is a small set $A \subseteq \mathfrak{C}$ of the monster model such that $\Th_A(\mathfrak{C})$ has the SB~property if and only if the theory $T$ is superstable and nonmultidimensional.

2. If $T$ has the SB~property and $A \subseteq \mathfrak{C}$ is any set, then $\Th_A(\mathfrak{C})$ also has the SB~property.
\end{conj}

Note that the second part of the last conjecture would follow immediately from Conjecture~\ref{main_SB_conj}.

Here is another very elementary-looking conjecture:

\begin{conj}
If $T$ does not have the SB~property, then there is an infinite collection of models of $T$ that are pairwise nonisomorphic and pairwise elementarily bi-embeddable.
\end{conj}

We are not sure why the last conjecture should be true, but a couple special cases of it have been verified below, namely the case where $T$ is unstable or has OTOP (Theorem~\ref{OP_NSB}) and the case where $T$ has a nomadic type (follows from Theorem~\ref{nomad_strong_NSB}).

\section{Some examples and counter-examples}

\begin{example}
The theory of all algebraically closed fields in the pure field language has the SB~property: any two fields that are bi-embeddable must have the same characteristic and the same absolute transcendence degree, and this is enough to determine an algebraically closed field up to isomorphism.  This is a natural example of an incomplete theory with the SB~property.
\end{example}

\begin{example}
A sort of generalization of the last example is the fact that any countable theory $T$ which is $\aleph_1$-categorical has the SB~property; this follows from Theorem~\ref{ttrans} above, or one can prove it directly using the classical analysis of such theories by Baldwin and Lachlan (\cite{bald_lach}).
\end{example}

\begin{example}
The theory of a dense linear ordering without endpoints does not have the SB~property; there are many ways to get counterexamples.  In fact, with a little thought one sees that \emph{any} theory of an infinite linear ordering, in the language with only ``$<$'', does not have the SB~property.
\end{example}

\begin{example}
\label{RCF}
The theory of real-closed fields, in the language $\left\{+, \cdot, <\right\}$, does not have the SB~property.  Given any model $M$ of this theory, let $\Inf(M)$ denote the positive non-standard elements of $M$, i.e. the elements of $M$ which are greater than every element of $\Q$.  Define the equivalence relation $\sim$ on $\Inf(M)$ by: $x \sim y$ if there is $n \in \N$ such that $x^n > y$ and $y^n > x$.  Note that $\sim$ partitions $\Inf(M)$ into convex sets.  By taking the real closures of the appropriate sets, we can find countable models $M$ and $N$ such that $\Inf(M) / \sim$ is order-isomorphic to $(\Q, <)$ and $\Inf(N) / \sim$ is order-isomorphic to the set of all rational points in $(0, 1]$.  Clearly $N$ and $M$ are non-isomorphic; to show that they are elementarily bi-embeddable, first note that $M$ and $N$ are isomorphic to subfields of one another, and then apply model completeness.
\end{example}

\begin{example}
\label{REF}
Let $T$ be the theory of infinitely many refining equivalence relations $E_i$ ($0 < i < \omega$) such that every $E_{i}$-class is split into two $E_{i+1}$-classes, and $E_i$ has exactly $2^i$ different classes.  Then $T$ eliminates quantifiers and is superstable and nonmultidimensional, but $T$ does not have the SB~property.  To prove this last statement, first note that the set $S^1(\acl^{eq}(\emptyset))$ of complete 1-types over $\acl^{eq}(\emptyset)$ naturally has the structure of an infinite binary tree whose nodes on the $i$th level correspond to the classes of $E_i$.  Then an arbitrary model $M \prec \mathfrak{C}$ of $T$ is completely determined by specifying, for each ``branch'' in $p \in S^1(\acl^{eq}(\emptyset))$, the cardinal number (possibly $0$) of realizations of $p$ in $M$.  The key point is that there exists an automorphism $f$ of $\mathfrak{C}$ whose orbits under its natural action on $S^1(\acl^{eq}(\emptyset))$ are infinite and dense in the Stone space topology: for each $i < \omega$, such an $f$ just needs to cyclically permute all the $E_i$-classes.  Now pick some arbitrary $p \in S^1(\acl^{eq}(\emptyset))$.  Let $M$ be a model such that $p, f(p), f^2(p), \ldots$ are each realized $2^{\aleph_0}$ times in $M$ and every other type in $S^1(\acl^{eq}(\emptyset))$ is omitted from $M$, and let $N$ be a model in which $p$ is realized countably many times, each of $f(p), f^2(p), \ldots$ is realized $2^{\aleph_0}$ times, and every other type from $S^1(\acl^{eq}(\emptyset))$ is omitted.  Then $M$ and $N$ are elementarily bi-embeddable (since $f(M) \prec N$ and $N \prec M$) but $M \ncong N$.
\end{example}

The next example shows that not all wandering types are nomadic:

\begin{example}
Let  $T = \Th(\Z; +)$.  $T$ is unidimensional, so it cannot have any nomadic types; but if $\widehat{\Z}$ is the profinite completion of $\Z$ (viewed as a ring), $\alpha, \beta \in \widehat{\Z}$ are algebraically independent, and $G \prec (\widehat{\Z}; +)$ is an elementary submodel such that $\alpha \cdot G = G$ and $\beta \notin G$, then $\tp(\beta / G)$ is a wandering type witnessed by the automorphism $f(x) = \alpha \cdot x$ of $G$.
\end{example}

The next example shows the difficulty of proving Conjecture~\ref{lifting-SB-conj}.  Even if $G$ is a weakly minimal group definable in a model $M$ and we have two models $H_1$ and $H_2$ of $\Th(G)$ and two elementary maps $f_1 : H_1 \rightarrow H_2$ and $f_2 : H_2 \rightarrow H_1$, there may not necessarily exist models $M_1$, $M_2$ such that $G(M_i) = H_i$ and elementary embeddings $f'_1 : M_1 \rightarrow M_2$ and $f'_2 : M_2 \rightarrow M_1$ that extend $f_1$ and $f_2$ respectively:

\begin{example}
\label{lifting_ex}
Let $G = \Z / {3 \Z} \oplus \left(\Z / {2 \Z} \right)^{\omega}$, and for each $i \in \omega$, let $G_i$ be the subgroup consisting of all elements whose $i$th coordinate in $\left(\Z / {2 \Z} \right)^{\omega}$ is zero.  Let $P = G \times G$, let $\pi : P \rightarrow G$ be the projection map onto the first coordinate, and let $f: G \times P \rightarrow P$ be the left regular action of $G$ on each fiber of $\pi$.  Let $P_i \subseteq P$ be the set of all $(g, h) \in P$ such that the $i$th projection of $h$ onto $\left(\Z / {2 \Z} \right)^{\omega}$ is zero.  Let $a \in \Z / {3 \Z} \oplus \prod_{i \in \omega} \Z / {2 \Z}$ be the element $(0, 1, 1, 1, \ldots)$ whose $\Z / {3 \Z}$ coordinate is $0$ and whose other coordinates are all $1$.  Finally, we define a map $h : P \rightarrow P$ as follows: if $(g,k) \in P$ and $g$ is divisible by $3$, then $h((g,k)) = (g,k)$; and if $(g,k) \in P$ and $g$ is not divisible by $3$, then $h((g,k)) = (-g, k + a)$.  Notice that $h^2$ is the identity map.  Let $T = \Th(G; +, G_i, \pi, f, h, P_i : i \in \omega)$.

Pick $H \models \Th(G)$ to be some model that omits $\tp(a)$, let $f_1 : H \rightarrow H$ be the identity, and let $f_2 : H \rightarrow H$ be the map $f_2(g) = -g$.  Then there are no models $M$ and $N$ of $T$ with elementary embeddings $f'_1 : M \rightarrow N$ and $f'_2 : N \rightarrow M$ such that $G(M) = H = G(N)$, $f'_1 \upharpoonright G(M) = f_1$, and $f'_2 \upharpoonright  G(N)  = f_2$.  (For if we had such models, suppose $g \in H$ is not divisible by $3$ and $(g, k) \in P(M)$.  If $f'_2 \circ f'_1 ((g,k)) = (-g, k')$ and $k'' \in G(M)$ is an element such that $f(k'', (-g, k')) = h((g, k))$, then $\tp(k'') = \tp(a)$, contradiction.)

\end{example}

\begin{example}
\label{single_const}
Here is an example of a theory $T$ which does not have the SB~property, but by adding a name for a single new element we get a theory with the SB~property.  $T$ is the theory of the structure $(\widehat{\Z_{(2)}}; S, \left\{E_i\right\}_{1 \leq i < \omega})$, where $\widehat{\Z_{(2)}}$ is the set of $2$-adic integers, $S$ is the unary successor function, and $(x, y) \in E_i$ if and only if $x$ and $y$ are congruent modulo $2^i$.  By the usual tricks, $T$ can be seen to eliminate quantifiers after adding a symbol for $0$.  From this one can show that $T$ is superstable and all types in $T$ have finite $U$-rank.  The presence of $S$ causes many types to be nonorthogonal, but there are still $2^{\aleph_0}$ nonorthogonality classes of minimal types.  Indeed, if $A$ is the set of all strong types over $\emptyset$ of single elements, then $A$ is homeomorphic to $\widehat{\Z_{(2)}}$, and if $p \in A$ corresponds to $a \in \widehat{\Z_{(2)}}$ and $x \in \widehat{\Z_{(2)}} - \Z$, then $a + x$ corresponds to a type in $A$ which is orthogonal to $p$.  Now we can use Theorem~\ref{nomadSB} below to show that $T$ does not have the SB~property.  However, if $a$ is any nonzero element of $\widehat{\Z_{(2)}}$, then $\acl^{eq}(a)$ contains names for every class of every $E_i$, and therefore any embedding between two models of $T$ that fixes $a$ will (when considered as an automorphism of the monster model) fix every nonorthogonality class of minimal types.  From this the SB~property of $\Th(\widehat{\Z_{(2)}}; S, \left\{E_i\right\}_{1 \leq i < \omega}, a)$ follows.
\end{example}

\section{The SB~Property and order properties}

In this section we prove that any theory in which a first-order formula can define arbitrarily long linear orderings does not have the SB~property.  In fact, we prove the more general result that any theory with an infinite linear order definable by a formula in $L_{\infty, \omega}$ (see Definition~\ref{inf-logics} below) does not have the SB~property.  In particular, a theory with the SB~property must be stable and have NOTOP.  We get our main result (Theorem~\ref{OP_NSB}) by modifying one of Shelah's many-model arguments (Theorem~VIII.3.2 of \cite{bible}), and we follow his proof closely.  The structure of the proof a little funny: instead of constructing two bi-embeddable models and showing that they are nonisomorphic, we construct a family of $2^\lambda$ bi-embeddable models for a suitably large infinite $\lambda$ and show that there must be a large subfamily of these models whose members are pairwise nonisomorphic.  In fact we do not know of any more direct proof of this theorem since we could not think of useful invariants of linear orderings that transfer to their corresponding Ehrenfeucht-Mostowski models; for instance, using the notation of Fact~\ref{EM_models}, it is possible that $I$ has uncountable cofinality while $EM(I)$ is an ordered structure with cofinality $\aleph_0$.

The models we build in this section are all EM models, and in addition we will use some infinitary combinatorics.  We emphasize that there are not any new deep ideas in this section and that we are only verifying that a previously-published construction works for our purposes.  All references in this section with roman numerals refer to results in \cite{bible}.  ``AP~n.m'' means ``Theorem/Claim/Lemma n.m of the Appendix of \cite{bible}.''

\begin{notation}
If $I$ is a linear order, then $I^*$ denotes the linear order which is the reverse of $I$, i.e. $\left\{ \langle b , a \rangle : \langle a , b \rangle \in I \right\}$.  If $I$ and $J$ are two linear orders defined on disjoint sets, then $I + J$ means the linear order in which each element of $I$ comes before every element of $J$, and $\sum_{i < \kappa} I_i$ has a similar meaning.
\end{notation}

\begin{definition}
\label{inf-logics}
1. For any language $L$ and any infinite cardinal $\kappa$, the logic $L_{\kappa, \omega}$ is the smallest set of formulas such that:

(a) $L_{\kappa, \omega}$ contains all first-order formulas in $L$;

(b) $L_{\kappa, \omega}$ is closed under the Boolean operations $\wedge$, $\vee$, and $\neg$ and universal or existential quantification over any variable;

(c) For any set $\left\{ \varphi_i(\bar{x}) : i \in I \right\}$ of formulas in $L_{\kappa, \omega}$ such that $\vert I \vert < \kappa$ and $\bar{x}$ has finite length, the conjunction $\bigwedge_{i \in I} \varphi_i(\bar{x})$ is in $L_{\kappa, \omega}$.

2. $L_{\infty, \omega}$ denotes $\bigcup_{\kappa \textup{ a cardinal}} L_{\kappa, \omega}$ (which is a proper class).
\end{definition}

\begin{definition}
\label{order-props}
1. (from \cite{unclass-NOTOP-unidim}) We say that a theory $T$ has the \emph{$L_{\infty, \omega}$ order property} if there is an $L_{\infty, \omega}$ formula $\psi(\bar{x}, \bar{y})$ (in the same vocabulary as $T$) such that for any linear ordering $\lambda$ there is a model $M$ of $T$ and a set $\langle \bar{a}_i : i < \lambda \rangle$ of tuples in $M$ such that $M \models \psi(\bar{a}_i, \bar{a}_j)$ if and only if $i < j$.

2. A theory with the $L_{\omega, \omega}$ order property is also called \emph{unstable}, and a theory that is not unstable is called \emph{stable}.
\end{definition}

A special case of the $L_{\infty, \omega}$ order property is:

\begin{definition}
(\cite{bible}, Definition~XII.4.2) $T$ has the \emph{omitting types order property} (or ``OTOP'' for short) if there is a type $p(\bar{x}, \bar{y}, \bar{z})$ such that for every ordinal $\lambda$ there is a model $M \models T$ and a sequence $\langle \bar{a}_i : i < \lambda \rangle$ from $M$ such that for any $i, j \in \lambda$, the type $p(\bar{a}_i, \bar{a}_j, \bar{z})$ is realized in $M$ if and only if $i < j$.  The negation of OTOP is called \emph{NOTOP}.
\end{definition}

\begin{remark}
Our definition of OTOP appears to be slightly weaker than the quoted definition in \cite{bible}, which says that the type $p(\bar{x}, \bar{y}, \bar{z})$ can be used to code \emph{any} binary relation in a model of $T$, but by the proof of Fact~\ref{EM_models} below these definitions are equivalent.
\end{remark}

\begin{lem}
\label{fields}
For any infinite regular cardinal $\kappa$, there exist linear orders $J_1$ and $J_2$ such that $\vert J_1 \vert = {\kappa}^+ = \vert J_2 \vert$, $\cf(J_1) = \kappa$, $\cf(J_2) = {\kappa}^+$, and $J_1$ and $J_2$ are bi-embeddable as orders.
\end{lem}

\begin{proof}
Let $T$ be the theory of real-closed fields in the language $\left\{0, +, \cdot,<\right\}$.  (The essential things are just that $T$ is some countable theory of ordered fields, $T$ is complete, and $T$ has definable Skolem functions.)  Let us recall the notation from Example~\ref{REF}: for $M \models T$, $\Inf(M)$ is the set of all positive elemetnts of $M$ which are greater than every element of $\Q$, and $x \sim y$ holds when there is $n \in \N$ such that $x^n > y$ and $y^n > x$.  By the usual compactness argument we can find a $<$-increasing indiscernible sequence $I$ in some model $M \models T$ such that $I \subseteq \Inf(M)$, no two distinct elements of $I$ are equivalent under $\sim$, and the order type of $I$ is ${\kappa}^+$.  Then $\dcl(I)$ is a model of $T$ of size ${\kappa}^+$, and if $J_2$ is the reduct of $\dcl(I)$ to the languange with only $<$ then $\cf(J_2) = {\kappa}^+$ (since $I$ is cofinal in $J_2$).  Let $J_1 \subseteq J_2$ be some sub-order of the form $J_1 = J' + J''$, where $J'$ is a nonempty bounded open interval in $J_1$ and $J''$ is some increasing sequence of elements above $J'$ such that $\cf(J'') = \kappa$.  Since $J_2$ is the order type of an ordered field, it is order-isomorphic to any of its nonempty open subintervals, in particular to $J'$, so $J_2$ is embeddable into $J_1$ and $\vert J_1 \vert = {\kappa}^+$.  Also clearly $\cf(J_1) = \kappa$, so we are done.
\end{proof}

We recall Definition~VIII.3.1:

\begin{definition}

1. If $(I, <)$ is an order and $J \subseteq I$, then for $\bar{a}, \bar{b} \in I$ we use the abbreviation ``$\bar{a} \sim \bar{b}$ mod $J$'' to mean that $\bar{a}$ and $\bar{b}$ have the same \emph{quantifier-free} types in the language of orders over the parameter set $J$.

2. If $(I, <)$ is an order, then an $I$-indexed sequence $\langle \bar{a}_s : s \in I \rangle$ of finite tuples from $M$ is \emph{skeletal} (in $M$) if it is indiscernible, and for any $\bar{c} \in M$, there is a finite $J \subseteq I$ such that if $\bar{s}, \bar{t} \in I$ and $\bar{s} \sim \bar{t}$ mod $J$, then for any $\varphi(\bar{x}, \bar{y})$ with the right numbers of variables, $M \models \varphi(\bar{a}_{\bar{s}}, \bar{c}) \leftrightarrow \varphi(\bar{a}_{\bar{t}}, \bar{c})$.

\end{definition}   

For example, the indiscernible sequences used to construct Ehrenfeucht-Mostowski models are skeletal in their Skolem hulls.

\begin{definition}
\label{contradictory}
(Definition~VIII.3.2 of \cite{bible}) The orders $I$ and $J$ are \emph{contradictory} if there is no order $K$ with infinite cofinality such there is a model $M$ with an antisymmetric definable relation $R$ on pairs of $m$-tuples and a sequence $\langle \bar{a}_s : s \in K + I^* + J^* \rangle$ of $m$-tuples from $M$ such that:

	1. $\langle \bar{a}_s : s \in K + I^* \rangle$ and $\langle \bar{a}_s : s \in K + J^* \rangle$ are both skeletal in $M$;
	
	2. If $s$ and $t$ are both in $K + I^*$ or both in $K + J^*$, then $M \models R(\bar{a}_s , \bar{a}_t)$ if and only if $s < t$.

\end{definition}

Note that in Definition~\ref{contradictory}, the indiscernibility of the sequences $\langle \bar{a}_s : s \in K + I^* \rangle$ and $\langle \bar{a}_s : s \in K+J^* \rangle$ implies that $R(\bar{x}, \bar{y})$ defines a linear ordering on each of these sequences.  Similarly the sequences of $\bar{a}_s$'s and $\bar{b}_s$'s in Definition~\ref{strongly_contradictory} below are each linearly ordered by $R$.

Definition~\ref{contradictory} is superficially quite similar to saying, ``any EM models built over indiscernibles with order types $I$ and $J$ must be nonisomorphic.''  However, there are examples of contradictory orders $I$ and $J$ and EM maps (see Fact~\ref{EM_models} below) such that $EM(I) \cong EM(J)$: for instance, if $I$ is an ordering of cofinality $\aleph_1$ and $J = I + \Q$, then $I$ and $J$ are contradictory by Lemmas~\ref{strong2weak} and \ref{distinct_cofs}, but it is possible that $EM(I) \cong EM(J)$ if there are constants in the Skolemized language that add a sequence of $\Q$-ordered elements above the skeleton of $EM(I)$.

\begin{definition}
\label{strongly_contradictory}
(Definition~AP~3.1 of \cite{bible}) The orders $I$ and $J$ are \emph{strongly contradictory} if they each have infinite cofinality and there are no orders $I'$ and $J'$ with infinite cofinality such that there is a model $M$ with an anti-symmetric definable relation $R$ on pairs of $m$-tuples and skeletal sequences $\langle \bar{a}_s : s \in I' + I^* \rangle$ and $\langle \bar{b}_s : s \in J' + J^* \rangle$ such that:

	1. For every $t \in J^*$ and $s' \in I'$, there is $s \in I'$ such that $s' < s$ and $M \models R(\bar{a}_s , \bar{b}_t)$,
	
	2. For every $t \in I^*$ and $s' \in J'$, there is $s \in J'$ such that $s' < s$ and $M \models R(\bar{b}_s , \bar{a}_t)$.
	
\end{definition}

\begin{remark}
It is unclear whether there are sequences that are contradictory but not strongly contradictory.  The more technical notion of being strongly contradictory seems to be needed only for the proof of Lemma~\ref{contradictory_sums} below.
\end{remark}

\begin{lem}
\label{strong2weak}
If $I$ and $J$ strongly contradictory, then they are contradictory.
\end{lem}

\begin{proof}

If $I$ and $J$ are not contradictory, as witnessed by the order $K$, the model $M$, and $\langle \bar{a}_s : K+ I^* + J^* \rangle$, then $M$ also witnesses the fact that $I$ and $J$ are not strongly contradictory, by letting (in the notation of Definition~\ref{strongly_contradictory}) $I' = K = J'$, using the same $\bar{a}_s$ for $s \in K + I^*$, and letting $\bar{b}_s = \bar{a}_s$ when $s \in K + J^*$.

\end{proof}

We now list a series of lemmas about contradictory orders from \cite{bible} that we need for our argument.

\begin{lem}
\label{distinct_cofs}
(AP~3.5) If $\cf(J) > \cf(I) \geq \aleph_0$, then $I$ and $J$ are strongly contradictory.
\end{lem}

\begin{lem}
\label{contradictory_sums}
(AP~3.6) Suppose that $\lambda$ is a regular uncountable cardinal, $\langle I_{\alpha} : \alpha < \lambda \rangle$ and $\langle J_{\alpha} : \alpha < \lambda \rangle$ are sequences of orderings such that each $I_\alpha$ and each $J_\alpha$ has infinite cofinality, $I = \sum_{\alpha < \lambda} I^*_{\alpha}$, $J = \sum_{\alpha < \lambda} J^*_{\alpha}$, and the set
\\ \\
	$S = \left\{\alpha < \lambda : \cf(\alpha) \geq \aleph_0 \textup{ and } I_{\alpha}, J_{\alpha} \textup{ are strongly contradictory} \right\}$
\\ \\
is a stationary subset of $\lambda$.  Then $I$ and $J$ are strongly contradictory.
\end{lem}

\begin{thm}
\label{stationary_partition}
(\cite{kunen}, Theorem~II.6.11) If $\lambda$ is a regular uncountable cardinal and $A \subseteq \lambda$ is stationary, then $A$ can be partitioned into $\lambda$ pairwise disjoint stationary subsets.
\end{thm}

\begin{lem}
\label{reglambda}
(Essentially AP~3.7) For any uncountable regular cardinal $\lambda$ there is a family of $2^{\lambda}$ pairwise strongly contradictory orders of cardinality $\lambda$ such that any two of them are bi-embeddable.
\end{lem}

\begin{proof}
$A = \left\{\alpha : \alpha < \lambda, \cf(\alpha) = \aleph_0 \right\}$ is a stationary subset of $\lambda$, so by Theorem~\ref{stationary_partition} there are subsets $A_i \subseteq A$ for each $i < \lambda$ which are stationary and pairwise disjoint.  Let $J_1$ and $J_2$ be two linear orders of size $\aleph_1$ as in the conclusion of Lemma~\ref{fields}, so that $\cf(J_1) = \aleph_0$ and $\cf(J_2) = \aleph_1$ and $J_1$ and $J_2$ are bi-embeddable.  By Lemma~\ref{distinct_cofs}, $J_1$ and $J_2$ are strongly contradictory.  For any $W \subseteq \lambda$ and $\alpha < \lambda$, let $I_{W, \alpha}$ be a copy of $J_1$ if $\alpha \in \bigcup_{\beta \in W} A_{\beta}$, and otherwise let it be a copy of $J_2$.  Let $I_W = \sum_{\alpha < \lambda} I^{*}_{W, \alpha}$.  If $W, U$ are subsets of $\lambda$ and $\gamma \in W - U$, then

$\left\{\alpha : \cf(\alpha) \geq \aleph_0, I_{W, \alpha} \textup{ and } I_{U, \alpha} \textup{ are strongly contradictory} \right\}$

contains $A_{\gamma}$.  So by Lemma~\ref{contradictory_sums}, $I_W$ and $I_U$ are strongly contradictory.  Also $I_W$ and $I_U$ are bi-embeddable, since they are sums each of whose components are bi-embeddable.
\end{proof}

\begin{lem}
\label{contraorders}
For any uncountable regular cardinal $\lambda$ there is a family of $2^{\lambda}$ pairwise strongly contradictory orders of cardinality $\lambda$ such that any two of them are bi-embeddable.
\end{lem}

\begin{proof}
The case where $\lambda$ is regular was taken care of above in Lemma~\ref{reglambda}, and the case where $\lambda$ is singular can be proved just like AP~3.8 (by requiring that the families $K_{\alpha}$ in the proof are not just pairwise strongly contradictory, but also pairwise bi-embeddable, using Lemma~\ref{reglambda} above).
\end{proof}

The next task is to define a class of $EM$ models which have indiscernible, linearly ordered skeletons.

From now on in this section, we consider a theory $T$ with the $L_{\infty, \omega}$ order property, as witnessed by the formula $\psi(\bar{x}, \bar{y})$, where both $\bar{x}$ and $\bar{y}$ have length $m$.  Let $L^1 \supset L$ be an expansion of the language by a single $2m$-ary relation $R(\bar{x}, \bar{y})$ and let $\textbf{K}^1$ be the class of all models of $T$ expanded to $L^1$ by interpreting $R(\bar{x}, \bar{y})$ by $\psi(\bar{x}, \bar{y})$.  If $T'$ is a theory in a language $L' \supseteq L$ and $\Gamma$ is a set of types in $T'$, then $PC_{L}(T', \Gamma)$ is the class of all reducts to $L$ of models of $T'$ that omit all types in $\Gamma$.

\begin{fact}
\label{Morleyization}
There is a first-order language $L^2 \supseteq L^1$, a theory $T^2$ in $L^2$, and a set $\Gamma$ of types in $T^2$ with the following properties:

1. $T^2 \supseteq \Th_{L^1}(\textbf{K}^1)$;

2. $T^2$ has Skolem functions;

3. For each subformula $\varphi(\bar{z})$ of $\psi$, there is a relation $R_{\varphi}(\bar{z}) \in L^2$;

4. $\textbf{K}^1 = PC_{L^1}(T^2, \Gamma)$;

5. For every model $M \in PC_{L^2}(T^2, \Gamma)$, every subformula $\varphi(\bar{z})$ of $\psi$, and every $\bar{a} \in M$, $M \models R_{\varphi}(\bar{a})$ if and only if $M \models \varphi(\bar{a})$.
\end{fact}

\begin{proof}
This follows from a standard technique of presenting the class of models of a sentence in $L_{\infty, \omega}$ as a pseudo-elementary class: $\Gamma$ is the set of all types of elements in nonstandard interpretations of the $R_{\varphi}$'s.  Also, $T^2$ will imply that $\forall \bar{x} \forall \bar{y} \left[R(\bar{x} , \bar{y}) \leftrightarrow R_{\psi}(\bar{x}, \bar{y}) \right]$.  For details, see the proof of Theorem~7.1.7 in \cite{baldwin_AEC}.
\end{proof}

\begin{fact}
\label{EM_models}
With the notation of Fact~\ref{Morleyization}, there is a mapping $EM^2$ from the class of all linear orderings into $PC_{L^2}(T^2, \Gamma)$ with the following properties:

	1. $\vert EM^2(I) \vert \leq \vert T^2 \vert + \vert I \vert$;

	2. For each linear order $I$, the model $EM^2(I)$ contains an indiscernible sequence $\langle \bar{a}_s : s \in I \rangle$ of $m$-tuples which is skeletal in $EM^2(I)$, called its \emph{skeleton};
	
	3. For each $I$, the skeleton $\langle \bar{a}_s : s \in I \rangle$ has the property that $EM^2(I) \models R(\bar{a}_s , \bar{a}_t)$ if and only if $s < t$;
	
	4. For any two orders $I$ and $J$, if $\langle \bar{a}_s : s \in I \rangle$ is the skeleton of $EM^2(I)$, $\langle \bar{b}_t : t \in J \rangle$ is the skeleton of $EM^2(J)$, and $\bar{s} \in I$ and $\bar{t} \in J$ are increasing finite tuples of the same length, then $\tp(\bar{a}_{\bar{s}}) = \tp(\bar{b}_{\bar{t}})$;
	
	5. If $J$ is a subordering of $I$, then there is an elementary embedding $f$ of $EM^2(J) \upharpoonright L$ into $EM^2(I) \upharpoonright L$.
	
We write $EM^1(I)$ for $EM^2(I) \upharpoonright L^1$ and $EM(I)$ for $EM^2(I) \upharpoonright L$.  Note that the skeleton of $EM^2(I)$ is still skeletal in $EM^1(I)$.
\end{fact}

\begin{proof}

This is standard so we do not go into too many details.  The main idea is that if $X = \langle \bar{a}_i : i \in I \rangle$ is a sufficiently long sequence in a model in $\textbf{K}^1$ which is $R$-increasing (that is, $i < j \Rightarrow R(\bar{a}_i, \bar{a}_j)$), then there is an infinite sequence $Y = \langle \bar{b}_i : i < \omega \rangle$ of $L^2$-indiscernibles in a model of $T^2$ such that the $L^1$-type of any finite increasing subsequence of $Y$ is realized by some finite $R$-increasing susequence of $X$.  ($X$ just has to be at least as long as the Hanf number $\beth_{\left(2^{\vert T^2 \vert}\right)^{+}}$; then use the Erd\H{o}s-Rado theorem and compactness to get $Y$.)  Now for any linear order $J$, we let $EM^2(J)$ be the Skolem hull of an $L^2$-indiscernible sequence with the same diagram as $Y$.  Full details are in sections VII.1 and VII.2 of \cite{bible} and in the proof of Theorem~10.11 in \cite{baldwin_AEC}.

\end{proof}

\begin{thm}
\label{OP_NSB}
If $T$ is complete and has the $L_{\infty, \omega}$ order property then $T$ does not have the SB~property.  In fact, for any cardinal $\kappa$, there is a cardinal $\lambda > \kappa$ such that $T$ has models $M_i$ for $i < 2^{\lambda}$, each having size $\lambda$, and for $i \neq j$, $M_i$ is not isomorphic to $M_j$ but $M_i$ is elementarily embeddable into $M_j$ (and vice-versa, by symmetry).
\end{thm}

\begin{proof}
We follow the argument for Theorem~VIII.3.2, making some modifications and clarifications, and proving the ``in fact...'' clause at the same time as the first clause.  Fortunately for what we are interested in we only have to deal with ``Case~I'' of Shelah's proof.

We let $\lambda = \left( {\kappa}^{+} + \vert T^2 \vert \right)^{\aleph_0}$.

\begin{claim}
\label{cardinal_computation}
(like AP~3.1 and 3.2(1))

1. There is no family $S$ of subsets of $\lambda$ such that $\vert S \vert = 2^{\lambda}$ and $S$ is a ``$(\lambda, \aleph_0)$-family,'' that is, each element of $S$ has size $\lambda$ and the intersection of any two distinct elements of $S$ has size less than $\aleph_0$.

2. Furthermore, if $2^{\lambda}$ is singular, there is a regular ordinal $\mu < 2^{\lambda}$ such that there is no $(\lambda, \aleph_0)$-family $S$ of subsets of $\lambda$ with $\vert S \vert = \mu$.
\end{claim}

\begin{proof}
1. Otherwise, if $S = \left\{A_i : i < 2^{\lambda} \right\}$, choose $B_i \subseteq A_i$ such that $\vert B_i \vert = \aleph_0$.  Then if $i \neq j$, $B_i \neq B_j$, so $\lambda^{\aleph_0} = \vert \left\{ B \subseteq \lambda : \vert B \vert = \aleph_0 \right\} \vert \geq \vert \left\{B_i : i < 2^{\lambda} \right\} \vert = 2^{\lambda}$.  But this contradicts the fact that $\lambda^{\aleph_0} = \lambda$.

2. Let $\mu_0 = \cf(2^{\lambda})$ and write $2^{\lambda} = \sum_{i < \mu_0} \mu_i$, where each $\mu_i$ is regular and less than $2^{\lambda}$.  Suppose 2 is false.  Then there is a $(\lambda, \aleph_0)$-family $S' = \left\{A_i : i < \mu_0 \right\}$ of subsets of $\lambda$.  Since there are bijections between each $A_i$ and $\lambda$, for each $i < \mu_0$ we can again choose $(\lambda, \aleph_0)$-families $S_i$ of $\mu_i$ subsets of $A_i$.  Thus the family $S = \bigcup_{i < \mu_0} S_i$ is a $(\lambda, \aleph_0)$-family of subsets of $\lambda$ and $\vert S \vert = 2^{\lambda}$, contradicting 1.
\end{proof}

Instead of using AP~3.3 as in VIII.3.2, we apply Lemma~\ref{contraorders} above to get pairwise contradictory orders $I_{\alpha}$ ($\alpha < 2^{\lambda}$), each of size $\lambda$, which are pairwise bi-embeddable.  For each $\alpha < 2^{\lambda}$ and $i \leq \lambda$, let $I^i_{\alpha}$ be an order isomorphic to $(I_{\alpha})^*$, and for future use we select elements $s^i_{\alpha} \in I^i_{\alpha}$.  Let $J_{\alpha} := \sum_{i \leq \lambda} I^i_{\alpha}$.  Then the $J_{\alpha}$'s are also pairwise bi-embeddable, and bi-embeddability is still preserved when we take the Skolem hulls $M_{\alpha} := EM^1(J_{\alpha})$ of these orders as in Fact~\ref{EM_models}.

\begin{claim}
If $\alpha, \beta < 2^{\lambda}$ and $M_\alpha \upharpoonright L \cong M_{\beta} \upharpoonright L$, then $M_{\alpha} \cong M_{\beta}$.
\end{claim}

\begin{proof}
Recall that $M_{\alpha}$ and $M_{\beta}$ are members of $PC_{L^1}(T^2, \Gamma)$, so by Fact~\ref{Morleyization}, the symbol ``$R$'' must have the same interpretation as the (infinitary) $L$-formula $\psi$ in both $M_{\alpha}$ and $M_{\beta}$.  So since the $L$-isomorphism between $M_{\alpha} \upharpoonright L$ and $M_{\beta} \upharpoonright L$ must respect all $L$-definable relations, the result follows.
\end{proof}

The last claim means that it suffices to find $2^\lambda$ of the $M_{\alpha}$'s which are nonisomorphic as structures in the expanded language $L^1$.  Suppose towards a contradiction that the number of nonisomorphic $M_{\alpha}$'s is less than $2^{\lambda}$, and let $\mu \leq 2^{\lambda}$ be regular.  Then for some model $M$, the set $S = \left\{ \alpha < 2^{\lambda} : M_{\alpha} \cong M \right\}$ has cardinality at least $\mu$, witnessed by isomorphisms $f_{\alpha}: M_{\alpha} \rightarrow M$.  For every $\alpha \in S$ let $\langle \bar{a}^{\alpha}_s : s \in J_{\alpha} \rangle$ be skeletal in $M_{\alpha}$.  Then for every $\alpha \in S$ the sequence $\langle \bar{b}^{\alpha}_s : s \in J_{\alpha} \rangle$ defined by $\bar{b}^{\alpha}_s = f_{\alpha}(\bar{a}^{\alpha}_s)$ is skeletal in $M$.  Also, for $\alpha \in S$, let $W_{\alpha} = \left\{\bar{b}^{\alpha}_s : s = s^i_{\alpha}, i < \lambda \right\}$.

\begin{claim}
\label{almost_disjointness}
If $\alpha, \beta$ are distinct elements of $S$ then $\vert W_{\alpha} \cap W_{\beta} \vert$ is finite.
\end{claim}

\begin{proof}
If not, then we can pick two strictly increasing functions $i, j : \omega \rightarrow \lambda$ such that for all $\xi < \omega$, if we let $s(\xi) = s^{i(\xi)}_{\alpha}$ and $t(\xi) = s^{j(\xi)}_{\beta}$, $\bar{b}^{\alpha}_{s(\xi)} = \bar{b}^{\beta}_{t(\xi)}$.  Let $\delta(1) = \sup \left\{i(\xi) : \xi < \omega \right\} \leq \lambda$ and $\delta(2) = \sup \left\{j(\xi) : \xi < \omega \right\} \leq \lambda$.  Now we prove a subclaim:

\begin{claim}
\label{skeletality}
The two sequences $$\langle \bar{b}^{\alpha}_s : s = s(\xi) \textup{ and } \xi < \omega, \textup{ or } s \in I^{\delta(1)}_{\alpha} \rangle$$ and $$\langle \bar{b}^{\beta}_s : s = t(\xi) \textup{ and } \xi < \omega, \textup{ or } s \in I^{\delta(2)}_{\beta} \rangle$$ are both skeletal in $M$.

\end{claim}

\begin{proof}

Let $\bar{c}$ be any finite tuple from $M$.  Since $\langle \bar{b}^{\alpha}_s : s \in J_{\alpha} \rangle$ is skeletal in $M$, there is a finite subset $K \subseteq J_{\alpha}$ such that if $\bar{s}, \bar{t} \in J_{\alpha}$ and $\bar{s} \sim \bar{t}$ mod $K$ then $(\bar{a}_{\bar{s}}, \bar{c})$ and $(\bar{a}_{\bar{t}}, \bar{c})$ have the same type in $M$.  Let $$K_1 = \left\{s(\xi) : \xi \in \omega \textup{ and there is an element of } K \textup{ between } s(\xi) \textup{ and } s(\xi + 1) \right\},$$ $$K_2 = \left\{s(\xi+1) : \xi \in \omega \textup{ and there is an element of } K \textup{ between } s(\xi) \textup{ and } s(\xi + 1) \right\},$$ and $K_3 = K \cap I^{\delta(1)}_{\alpha}.$  Then $K' = K_1 \cup K_2 \cup K_3$ is finite.  If $\bar{s}, \bar{t}$ are finite sequences from $\ran(s) \cup I^{\delta(1)}_{\alpha}$ and $\bar{s} \sim \bar{t}$ mod $K'$, then $\bar{s} \sim \bar{t}$ mod $K$, so $\tp(\bar{a}_{\bar{s}}, \bar{c}) = \tp(\bar{a}_{\bar{s}}, \bar{c})$, proving that $K'$ witnesses the skeletality of the first sequence.  The proof that the second sequence is skeletal is very similar.

\end{proof}

But Claim~\ref{skeletality} implies that $(I^{\delta(1)}_{\alpha})^* \cong I_{\alpha}$ and $(I^{\delta(2)}_{\beta})^* \cong I_{\beta}$ are not contradictory, which is a contradiction, finishing the proof of Claim~\ref{almost_disjointness}.

\end{proof}

By Claim~\ref{almost_disjointness}, $\left\{W_{\alpha} : \alpha \in S \right\}$ is a $(\lambda, \aleph_0)$-family of subsets of $M^m$, which has size $\lambda$, and $\vert S \vert = \mu$.  Since the argument for Claim~\ref{almost_disjointness} did not depend at all on the choice of the regular $\mu \leq 2^{\lambda}$, this contradicts either part~1 of Claim~\ref{cardinal_computation} if $2^{\lambda}$ is regular, or else part~2 of Claim~\ref{cardinal_computation} if $2^{\lambda}$ is singular.

\end{proof}

\begin{quest}
If $T$ is unstable, then can we strengthen the conclusion of Theorem~\ref{OP_NSB} by adding the clause: ``each $M_i$ is $\kappa$-saturated''?  In Theorem~VIII.3.2 of \cite{bible}, Shelah does construct his models $M_i$ to be $\kappa$-saturated, but he accomplishes this with sequences that are ``$\F^b_{\kappa}$-constructible,'' and this particular notion of isolation does not admit saturated models in the sense of Definition~IV.1.1 of \cite{bible}.  What we want is some notion of isolation $\F_{\kappa}$ such that even in an unstable theory $\F_{\kappa}$ satisfies enough of the axioms of isolation (as in section~IV.1 of \cite{bible}) so that over any set there is a model which is both $\F_{\kappa}$-constructible and $\F_{\kappa}$-saturated.  Furthermore, $\F_{\kappa}$-saturated models should also be $\kappa$-saturated in the usual sense.
\end{quest}

\begin{remark}
The reason why we consider only the $L_{\infty, \omega}$ order property in Theorem~\ref{OP_NSB}, instead of the more general class of $L_{\infty, \infty}$-definable orders, is that we do not know whether this more general order property would imply the existence of a functor from linear orders into models as in Fact~\ref{EM_models}.  If we had such a functor, the rest of the argument for Theorem~\ref{OP_NSB} would go through.
\end{remark}

\section{The SB~property in the stable case}

\subsection{The strictly stable case}

In this subsection, we assume that $T$ is a stable theory (see Definition~\ref{order-props}), and aim towards proving that the SB~property cannot hold for \emph{strictly stable} theories, that is:

\begin{definition}
1. The theory $T$ is \emph{strictly stable} if $T$ is stable and there is an infinite chain $p_0(\bar{x}) \subseteq p_1(\bar{x}) \subseteq \ldots$ of types such that for every $i < \omega$, $p_{i+1}(\bar{x})$ forks over $\dom(p_i)$.

2. The theory $T$ is \emph{superstable} if it is stable but not strictly~stable.
\end{definition}

\begin{definition}
\label{f-constr}
(From \cite{bible}, section~IV.2) 

1. A type $p \in S(A)$ is \emph{$f$-isolated} if there is a finite set $B \subseteq A$ such that $p$ does not fork over $B$.  (In \cite{bible} this is called ``$\F^f_{\aleph_0}$-isolation.'')

2. A set $B$ is \emph{$f$-constructible over $A$} if $A \subseteq B$ and there is an enumeration $\left\{b_i : i < \alpha \right\}$ of all the elements of $B - A$ such that for each $i < \alpha$, $\tp(b_i / A \cup \left\{b_j : j < i \right\})$ is $f$-isloated.

\end{definition}

The next two lemmas about $f$-constructible models were proved already in \cite{bible} in a much more general setting, but since the proofs are so short it will not hurt to spell them out here for ease of reference.

\begin{lem}
\label{f-constructions}
For any set $A$, there is a model $B \supseteq A$ such that $B$ is $f$-constructible over $A$ and $\vert B \vert \leq \vert A \vert + \vert T \vert$.
\end{lem}

\begin{proof}
First, enumerate all consistent formulas $\varphi(\bar{x}, \bar{a})$ over $A$ and realize each one successively by some element whose type over $A$ and all the previously chosen elements does not fork over $\bar{a}$.  The union of this sequence is a set $A_1 \supseteq A$ such that $A_1$ is $f$-constructible over $A$ and $\vert A_1 \vert \leq \vert A \vert + \vert T \vert$.  Similarly we can construct sets $A_2 \subseteq A_3 \subseteq \ldots$ such that $A_{i+1}$ realizes every consistent formula over $A_i$.  $\bigcup_{i<\omega} A_i$ is the model we want.
\end{proof}

Note that in a stable theory, $f$-isolation satisfies most of the usual properties of a notion of isolation (symmetry, transitivity, etc.), except that there is no such thing as an ``$f$-saturated model''  and the model constructed in the previous lemma is not prime in any reasonable class of models.

\begin{lem}
\label{f-atomicity}
If the model $M$ is $f$-constructible over $A$, then for any $\bar{b} \in M$, there is a finite set $F \subseteq A$ such that $\tp(\bar{b}/A)$ does not fork over $F$.
\end{lem}

\begin{proof}
Label the $f$-construction of $M$ as $\left\{b_i : i < \alpha \right\}$ as in Definition~\ref{f-constr}, and let $\bar{b} = (b_{i_1}, \ldots, b_{i_n})$ for some ordinals $i_1 < \ldots < i_n < \alpha$.  The proof goes by induction on $i_n$, so assume that the lemma has been proved for all sequences from $M$ that are constructed before the $i_n$-th stage.  By the definition of $f$-constructability, there are $j_1 < \ldots < j_m < i_n$ and a finite $A_0 \subseteq A$ such that $\tp(b_{i_n} / A \cup \left\{b_k : k < i_n\right\})$ does not fork over $A_0 \cup \left\{b_{j_i}, \ldots, b_{j_m}\right\}$.  Let $b' = (b_{j_1}, \ldots, b_{j_m}, b_{i_1}, \ldots, b_{i_{n-1}})$.  By our induction hypothesis, there is also a finite $A_1 \subseteq A$ such that $\tp(b'/ A)$ does not fork over $A_1$.  Let $F = A_0 \cup A_1$.  Then by monotonicity and symmetry, $A \ind_F b'$ and $A \ind_{F b'} b_{i_n}$, so by transitivity, $A \ind_F b' b_{i_n}$.  By symmetry and monotonicity again, we get the conclusion we want.
\end{proof}

Leo Harrington deserves thanks for pointing out a way to both simplify and generalize my original argument for the next theorem.

\begin{thm}
\label{nonsstable_SB}
If $T$ is strictly stable then $T$ does not have the SB~property.
\end{thm}

\begin{proof}
By the hypothesis there is an $\omega$-indexed sequence of types $p_i(\bar{x}, \bar{a}_i)$ over parameters $\bar{a}_i$, with $\bar{a}_i$ a possibly infinite tuple but $\bar{x}$ finite, such that for every $i < \omega$, the following hold:

1. $\bar{a}_i \subsetneq \bar{a}_{i+1}$;

2. $p_i(\bar{x}; \bar{a}_i) \subseteq p_{i+1}(\bar{x}; \bar{a}_{i+1})$;

3. $p_{i+1}(\bar{x}; \bar{a}_{i+1})$ forks over $\bar{a}_i$.

By passing to nonforking extensions, we may also assume that every type $p_i(\bar{x}, \bar{a}_i)$ is complete.  Let $A = \bigcup_{i < \omega} \bar{a}_i$ and let $p = \bigcup_{i < \omega} p_i(\bar{x}; \bar{a}_i)$ (so $p \in S(A)$).  By stability, $T$ has a saturated model $M$ of size at least $\vert T \vert^{+}$, and we can pick $M$ so that $\tp(A/M)$ does not fork over $\emptyset$.  By Lemma~\ref{f-constructions}, there is a model $N$ such that $\vert N \vert = \vert M \vert$ and $N$ is $f$-constructible over $M \cup A$.

We claim that the model $N$ contains no realization of the type $p$.  Suppose towards a contradiction that $N$ contains $\bar{b}$ realizing $p$.  By Lemma~\ref{f-atomicity} there is some $i < \omega$ such that $\tp(\bar{b} / M \cup A)$ does not fork over $M \cup \bar{a}_i$.  From $A \ind M$ and monotonicity, we get $A \ind_{\bar{a}_i} M$, so by transitivity it follows that $A \ind_{\bar{a}_i} \bar{b} M$.  By symmetry and monotonicity, $\bar{b} \ind_{\bar{a}_i} A$.  But this contradicts the fact that $\bar{b}$ is supposed to realize the type $p \in S(A)$ which forks over $\bar{a}_i$.

Since $N$ does not realize the type $p$ and $\vert \dom(p) \vert  = \aleph_0 < \vert N \vert$, $N$ is not saturated and therefore it cannot be isomorphic to $M$.  But $M \prec N$ and $N$ is elementarily embeddable into $M$ because $M$ is saturated.

\end{proof}

\subsection{Nomadic types and the SB~property}

In this subsection, any theory $T$ mentionned is assumed to be stable, but not necessarily superstable unless specified.  We follow the convention that ``stationary types'' are really parallelism classes of types over small sets (see the next definition below).

The standard notions below are from \cite{bible}; for a more elementary exposition of them, see \cite{baldwin} or \cite{poizat}.

\begin{definition}

\begin{enumerate}

\item A type $p$, not necessarily complete, is \emph{stationary} if for every set $A$ there is only one nonforking extension of $p$ to $S(\dom(p) \cup A)$.

\item If $p$ is a stationary type then $p \vert A$ denotes the natural restriction to $S(A)$ of the nonforking extension of $p$ to $S(\dom(p) \cup A)$.  Usually we will only use the notation ``$p \vert A$'' when this resulting type is also stationary (in which case we say that $p$ is \emph{based on} $A$).

\item Two stationary types $p$ and $q$ are \emph{parallel} if $$p \vert \left(\dom(p) \cup \dom(q)\right) = q \vert \left(\dom(p) \cup \dom(q)\right).$$

\item If $p \vert A$ is stationary and $I$ is an independent set over $A$ such that each element of $I$ realizes $p \vert A$, then we say that $I$ is a \emph{Morley sequence based on $p \vert A$}.  (Note that such a sequence must be indiscernible over $A$.)

\item Two stationary types $p$ and $q$ are \emph{orthogonal}, written in symbols as ``$p \bot q$,'' if for any set $B \supseteq \dom(p) \cup \dom(q)$ and any tupes $\bar{a}$ realizing $p \vert B$ and $\bar{b}$ realizing $q \vert B$, $\tp(\bar{a} / B \cup \bar{b})$ does not fork over $B$.

\item A type $p$ is \emph{regular} if it is nonalgebraic and orthogonal to all of its forking extensions.
\end{enumerate}
\end{definition}

The next notion does not seem to have been discussed in the literature before, but we think that it is closely tied to the SB~property:

\begin{definition}
\label{nomads}
If $p$ is a stationary type in a stable theory, then $p$ is \emph{nomadic} if there exists an automorphism $f$ of $\mathfrak{C}$ such that for every $i < j < \omega$, $f^i(p) \bot f^j(p)$.  (Here we write ``$f^0$'' for the identity map.)
\end{definition}

\begin{remark}
For a type $p$ to be nomadic, it is not enough for $p$ to have infinitely many pairwise-orthogonal conjugates, even if $p$ is regular.  For example, let $T = \Th(2^{\omega}, \left\{E_i\right\}_{i < \omega})$, where for any two binary sequences $\sigma, \tau \in 2^{\omega}$, $\sigma E_i \tau$ holds if and only if $\sigma(i) = \tau(i)$.  Then any nonalgebraic stationary $1$-type has $2^{\aleph_0}$ pairwise-orthogonal conjugates, but for any automorphism $f$ of $\mathfrak{C}$, $f^2$ preserves all nonorthogonality classes.  Also, for the record, this theory does have the SB~property.
\end{remark}

The next lemma was extracted from the proof of Lemma~20.19 in \cite{poizat} (and a similar argument appeared earlier in \cite{bible}).

\begin{lem}
\label{stp-nomads}
If $p$ is a stationary regular type based on $A$ and there is an automorphism $g$ of the monster model such that $g(p) \bot p$ and $\stp(g(A) / \emptyset) = \stp(A / \emptyset)$, then $p$ is nomadic.
\end{lem}

\begin{proof}
Let $A_1, A_2, \ldots$ be a Morley sequence in the nonforking extension of $\stp(A / \emptyset)$ to $A \cup g(A)$.  Then $A_0 = A, A_1, A_2, \ldots$ and $g(A), A_1, A_2, \ldots$ are both Morley sequences in $\stp(A / \emptyset)$.  Let $p_i$ be the type over $A_i$ which corresponds to $p \in S(A)$ by an elementary map from $A$ onto $A_i$.  Since nonorthonality is a transitive relation on regular types and $p \bot g(p)$, it follows that $p_i \bot p$, and so by indiscernibility $p_i \bot p_j$ for any $i < j < \omega$.  Therefore $p$ is a regular nomad, as witnessed by an automorphism $f$ which sends each $A_i$ to $A_{i+1}$.
\end{proof}

\begin{cor}
\label{automorphnomads}
If $p$ is a stationary regular type based on $A$ and there are automorphisms $g$ and $g'$ of $\mathfrak{C}$ such that $g(p) \bot g'(p)$ and $\stp(g(A) / \emptyset) = \stp(g'(A) / \emptyset)$, then $p$ is nomadic.
\end{cor}

\begin{proof}
Since orthogonality is preserved by the automorphism $g^{-1}$, $p \bot \left[g^{-1} \circ g'\right](p)$.  Since $g^{-1} \circ g'$ fixes $\stp(A/\emptyset)$, the conclusion of the proposition follows from Lemma~\ref{stp-nomads}.
\end{proof}

\begin{prop}
\label{regnomads}
If $T$ is superstable and there is a nomadic type in $T$, then there is a regular nomadic type in $T$.
\end{prop}

\begin{proof}
Suppose $p$ is a nomadic type, as witnessed by a map $f \in \Aut(\mathfrak{C})$ such that $p \bot f(p) \bot \ldots$.  It is folklore that superstability implies that $p$ is domination equivalent to a finite product $r_1 \otimes \ldots \otimes r_n$ of regular types -- see, for instance, Theorem XIII.3.6 of \cite{baldwin} or Corollary~1.4.5.7 of \cite{pillay}.  It follows immediately that for any $i$, $f^i(p)$ is domination equivalent to $f^i(r_1 \otimes \ldots \otimes r_n)$.  Since orthogonality is invariant under the substitution of domination-equivalent types, it follows that for any $i \neq j$, $$\otimes_k f^i(r_k) = f^i(\otimes_k r_k)  \bot f^j(\otimes_k r_k) = \otimes_k f^j(r_k).$$  Therefore, $f^i(r_1) \bot f^j(r_1)$, so $r_1$ is nomadic.
\end{proof}

Recall that a superstable theory is called \emph{multidimensional} (see \cite{bible}, V.5) if there is a set $X \subseteq S(\mathfrak{C})$ such that $\vert X \vert > 2^{\vert T \vert}$ and every pair of types in $X$ is orthogonal.  This notion has also been called ``nondimensionality'' \cite{poizat} and ``unboundedness'' \cite{buechler}.  It is possible for $T$ to be superstable and for $T$ to have a nomadic type without $T$ being multidimensional.  For example, let $T$ be as in Example~\ref{REF}.  Then $T$ can have no more than $2^{\aleph_0}$ pairwise-orthogonal types over any model; but, as we observed in Example~\ref{REF}, any nontrivial strong type over $\emptyset$ is nomadic.

\begin{lem}
\label{multidimnomad}
If $T$ is a superstable multidimensional theory then there is a regular nomadic type in $T$.
\end{lem}

\begin{proof}
By stability and multidimensionality, there are models of $T$ over which there are arbitrarily large collections of stationary and pairwise-orthogonal types.  Every stationary type has a base of size $< \vert T \vert^{+}$, so by the pigeonhole principle there must be a pair of orthogonal stationary types $p \in S(A)$ and $q \in S(B)$ over small sets $A$ and $B$ in the monster model such that $\textup{stp}(A / \emptyset) = \textup{stp}(B / \emptyset)$ and an automorphism $g$ of the monster model which preserves $\textup{stp}(A / \emptyset)$ and maps $p$ onto $q$.  The rest of the proof is the same as in Proposition~\ref{regnomads}: by superstability, $p$ is domination-equivalent to some finite product $r_1 \otimes \ldots \otimes r_n$ of regular types, and so $r_1 \otimes \ldots \otimes r_n \bot g(r_1 \otimes \ldots \otimes r_n) = g(r_1) \otimes \ldots \otimes g(r_n)$.  This implies that $r_1 \bot g(r_1)$, and Lemma~\ref{stp-nomads} finishes the proof.

\end{proof}

\begin{definition}
\label{dimension}
(\cite{bible}, Definitions III.1.6 and III.3.3) 
1. If $I$ and $J$ are two infinite indiscernible sequences, we say $I$ and $J$ are \emph{equivalent}, and write $I \sim J$, if the average types $\textup{Av}(I /  I \cup J)$ and $\textup{Av}(J / I \cup J)$ are equal.

2. If $I$ is any infinite indiscernible sequence and $A$ and $B$ are any two sets, then $\textup{dim}(I, A , B)$, the \emph{dimension of $I$ over $A$ in $B$}, is defined to be the minimum of $$\left\{\vert J \vert : J \sim I \textnormal{ and } J \textnormal{ is a maximal indiscernible sequence in } B \textnormal{ over } A\right\}.$$

3. We write $\textup{dim}(I, B)$ for $\textup{dim}(I, \emptyset, B)$.

4. If $p \in S(A)$ is a stationary type, then ``$\dim(p, A, B)$'' means $\dim(I, A, B)$ where $I$ is any Morley sequence based on $p$.
\end{definition}

\begin{fact}
\label{trueness}
(\cite{bible}, Lemma~III.3.9) If $M \models T$ and $\dim(I, A, M) \geq \kappa(T)$, then for \emph{any} sequence $J \subseteq M$ which is indiscernible over $A$ and equivalent to $I$, $\vert J \vert \leq \dim(I, A, M)$.
\end{fact}

\begin{remark}
Usually in Definition~\ref{dimension} we care about the case where $B$ is a model.  If $T$ is superstable, then any two maximal indiscernible sequences over $A$ in $B$ equvialent to $I$ will have the same cardinality (by Fact~\ref{trueness}).
\end{remark}

\begin{fact}
\label{basechange}
(\cite{bible}, Corollary~III.3.5) If $I$ is an indiscernible set over $A$, $B$ is any other set, and $\vert I \vert > \kappa(T) + \vert B \vert$, then there is a subset $J \subseteq I$ such that $\vert J \vert \leq \kappa(T) + \vert B \vert$ and $I - J$ is indiscernible over $A \cup B \cup J$.
\end{fact}

\begin{lem}
\label{basechange2}
If $A \subseteq B \subseteq C$ then $\dim(I , B, C) \leq \dim(I, A , C)$.  If $A \subseteq B \subseteq M$, $M \models T$, and $\dim(I, B, M) > \kappa(T) + \vert B - A \vert$, then $\dim(I, A, M) = \dim(I, B, M)$.
\end{lem}

\begin{proof}
The first statement follows immediately from the definitions.  For the second statement, suppose that there was an $A$-indiscernible sequence $J \subseteq M$ such that $J$ is equivalent to $I$ and $\vert J \vert > \dim(I, B, M)$.  Then by Fact~\ref{basechange}, there is a subsequence $J' \subseteq J$ such that $\vert J' \vert = \vert J \vert$ and $J'$ is indiscernible over $B$.  But by Fact~\ref{trueness}, $\vert J' \vert \leq \dim(I, B, M)$, contradiction.
\end{proof}

\begin{lem}
\label{morley_basechange}
If $I$ is an infinite Morley sequence over the set $A$ and $B$ is a set such that $\vert B \vert + \kappa_r(T) < \vert I \vert$, then there is a subsequence $I' \subseteq I$ such that $\vert I' \vert < \vert B \vert + \kappa_r(T)$ and $(I - I') \ind_{A} A B I'$, and thus $I - I'$ is a Morley sequence over $A \cup B \cup I'$.
\end{lem}

\begin{proof}
Pick $I_0 \subseteq I$ such that $\vert I_0 \vert \leq \vert B \vert + \kappa_r (T)$ and $B \ind_{A I_0} A I$.  Then symmetry and monotonicity of nonforking imply that $(I - I_0) \ind_{A I_0} B$, so by the fact that $(I - I_0) \ind_A I_0$ plus transitivity, $(I - I_0) \ind_{A} A B$.  Repeating the same argument $\omega$ times, we get a series of subsequences $I_i$, with $I_{i+1} \subseteq I - \bigcup_{j \leq i} I_j$, such that $I - (I_0 \cup \ldots \cup I_{i+1})$ is free from $B \cup \bigcup_{j \leq i} I_j$ over $B \cup \bigcup_{j < i} I_j$.  By transitivity of nonforking, $I' = \bigcup_{i < \omega} I_i$ satisfies the conclusion of the lemma.
\end{proof}

To build counterexamples to the SB~property, we will make use of the following notion of isolation:

\begin{definition}
1. (\cite{bible}, Definition~IV.2.1) If $\kappa$ is an infinite cardinal, the type $p \in S(A)$ is \emph{$\F^a_{\kappa}$-isolated} if there is a set $B \subseteq A$ such that $\vert B \vert < \kappa$ and there is a type $q \in S(\acl^{eq}(B))$ such that $q \vdash p$.

2. A model $M$ is called \emph{$\F^a_{\kappa}$-saturated} if for every $B \subseteq M$ such that $\vert B \vert < \kappa$ and every type $q \in S(\acl^{eq}(B))$, $q$ is realized in $M$.

3. A model $M$ is called \emph{$\F^a_{\kappa}$-prime over $A$} if $A \subseteq M$, $M$ is $\F^a_{\kappa}$-saturated, and for any other $\F^a_{\kappa}$-saturated model $N$ containing $A$, there is an elementary map $f$ that fixes $A$ pointwise and sends $M$ into $N$.
\end{definition}

\begin{remark}
The most familiar example of $\F^a_{\kappa}$-isolation is when $\kappa = \kappa_r(T)$, in which case this is more commonly called ``$a$-isloation'' and $\F^a_{\kappa}$-saturated models are called ``$a$-models'' (see section 1.4.2 of \cite{pillay}).  We use the cardinal parameter below in order to get counterexamples to the SB~property that are $\kappa$-saturated for arbitrarily large $\kappa$, but the reader who does not care about his extra generality can ignore all the $\kappa$'s in the rest of this section and pretend we are only using $a$-isolation.
\end{remark}

\begin{fact} (IV.2 and IV.3 of \cite{bible})
Suppose $\kappa \geq \kappa(T)$.  Then over any set $A$ there exists a model that is $\F^a_{\kappa}$-prime over $A$.
\end{fact}

The next two technical-looking lemmas are crucial for our computations of dimensions of sequences inside saturated models.

\begin{fact}
\label{prime_dimensions}
(Theorem~IV.4.9(5) of \cite{bible})  If $\kappa \geq \kappa(T)$, $\kappa$ is regular, and $M$ is $\F^a_{\kappa}$-prime over $A$, then for any infinite indiscernible sequence $I$ over $A$ such that $I \subseteq M$, $\dim(I, A, M) \leq \kappa$.
\end{fact}

\begin{fact}
\label{regular_dimensions}
(Lemma~V.2.4(1) of \cite{bible}) If $\Av(I, I)$ and $\Av(J, J)$ are nonorthogonal regular types, $M$ is an $\F^a_{\kappa_r(T)}$-saturated model, and $I, J \subseteq M$, then $\dim(I, M) = \dim(J, M)$.
\end{fact}

The next theorem is really just a generalization of the idea discussed in Example~\ref{REF}.

\begin{thm}
\label{nomadSB}
If $T$ is a stable theory with a regular nomadic type then $T$ does not have the SB~property.  Moreover, for any cardinal $\kappa \geq \vert T \vert^{+}$, there are witnesses $M, N$ to the failure of the SB~property such that both $M$ and $N$ are $\kappa$-saturated.
\end{thm}

\begin{proof}
We will prove the ``moreover'' clause at the same time as the main statement.  Fix some regular $\kappa \geq \vert T \vert^{+}$ and pick a regular nomad $p$, with $f$ as in Definition~\ref{nomads}.  Pick a base $A$ for $p$ such that $\vert A \vert < \vert T \vert^{+}$.  To simplify notation, let $p_i = f^i(p)$ (which, note, is based on $f^i(A)$).  Let $B$ (for ``Base'') be the set $\bigcup_{i < \omega} f^i(A)$.  Pick some ordinal $\alpha$ such that $\aleph_{\alpha} > \kappa$, and for every $i < \omega$, let $I_i$ be a Morley sequence in $p_i \vert \left(B \cup \bigcup_{j < i} I_j\right)$ of length $\aleph_{\alpha + i}$.  Choose models $M$ and $N$ such that $M$ is $\textup{\textbf{F}}^a_{\kappa}$-prime over $B \cup \bigcup_{i < \omega} I_i$ and $N$ is $\textup{\textbf{F}}^a_{\kappa}$-prime over $B \cup \bigcup_{0< i < \omega} I_i$.  Note that there is an elementary map $g_1$ with domain $B \cup \bigcup_{i < \omega} I_i$ such that $g \upharpoonright B = f \upharpoonright B$ and $G$ maps each $I_i$ into $I_{i+1}$ (simply by the stationarity of $p \vert B$), and the inclusion map $g_2 : \left(B \cup \bigcup_{0 < i < \omega} I_i \right) \rightarrow \left(B \cup \bigcup_{i < \omega} I_i\right)$ is also elementary.  These maps can be extended to elementary embeddings $\overline{g_1} : M \hookrightarrow N$ and $\overline{g_2} : N \hookrightarrow M$ (by primality), so $M$ and $N$ are elementarily bi-embeddable.

Now suppose that $A' \subseteq N$ realizes $\textup{tp}(A / \emptyset)$ and $f'$ is an automorphism of the monster model mapping $A$ onto $A'$.  Let $p' = f'(p \vert A)$, and note that by the saturation of $N$ there is an infinite Morley sequence $I \subseteq N$ over $A'$ in the type $p'$.

\begin{claim}
If $I' \subseteq N$ is an indiscernible sequence based on $p'$, then $\textup{dim}(I', N) \neq \aleph_{\alpha}$.
\end{claim}

\begin{proof}

Case 1: For every $i$ such that $0 < i < \omega$, $p' \bot p_i$.  Suppose towards a contradcition that there is $J \subseteq N$ such that $\textup{Av}(J / J)$ is parallel to $p'$ and $\vert J \vert \geq \aleph_{\alpha}$.  Note that for any countable $J' \subseteq J$, $J - J'$ is a Morley sequence over $ J'$, parallel to $p'$.  By Lemma~\ref{morley_basechange}, there is a subset $J'' \subseteq (J - J')$ such that $J - (J' \cup J'')$ is free from $A' \cup B$ over $ J'$ and $\vert J'' \vert \leq \kappa_r(T) + \vert A' \cup B \vert \leq \vert T \vert^{+} < \aleph_{\alpha}$.  Let $J^* = J - (J' \cup J'')$.  Then $J^*$ is a Morley sequence over $A' \cup B \cup J'$ in the type $p' \vert (A' \cup B \cup J')$, and $\vert J^* \vert = \vert J \vert \geq \aleph_{\alpha}$.

For any nonzero $i \in \omega$, we can use Lemma~\ref{morley_basechange} again to get a subset $I'_i \subseteq I_i$ such that $\vert I'_i \vert \leq \kappa_r(T) + \vert A' \cup B \vert \leq \vert T \vert^+ < \aleph_{\alpha}$ and $(I_i - I'_i) \ind_{B} A' B I'_i$.  A third use of Lemma~\ref{morley_basechange} yields a set $\widehat{J} \subseteq J^*$ such that $\vert \widehat{J} \vert < \aleph_{\alpha}$ and $J^* - \widehat{J}$ does not fork with $\bigcup_{0 < i < \omega} I'_i$ over $A' \cup B \cup J'$.  Now since $p'$, the type on which $J^* - \widehat{J}$ is based, is orthogonal to every $p_i$ and the $p_i$'s are pairwise orthogonal, $J^* - \widehat{J}$ is free from $\bigcup_{0 < i < \omega} (I_i - I'_i)$ over $A' \cup B \cup J' \cup \bigcup_{0 < i < \omega} I'_i$.  Transitivity then implies that $J^* - \widehat{J}$ is a Morley sequence over $A' \cup B \cup \bigcup_{0 < i < \omega} I_i$.  But $\vert J^* - \widehat{J} \vert \geq \aleph_{\alpha} > \kappa$, so this contradicts Fact~\ref{prime_dimensions}.  So when $f'(p)$ is orthogonal to every $p_i$ where $i$ is nonzero, $\textup{dim}(I, N) < \aleph_{\alpha}$. 

Case 2: There is some nonzero $i < \omega$ such that $f'(p)$ is nonorthogonal to $p_i$.  Then by Fact~\ref{regular_dimensions}, $\textup{dim}(I, N) = \textup{dim}(J, N)$ for any indiscernible sequence $J \subseteq N$ whose average type is parallel to $p_i$.  By construction, $\dim(J, N) \geq \dim(J, f_i(A), N) > \aleph_{\alpha}$, so the claim is proved.

\end{proof}

Now in $M$ there is an indiscernible sequence $I$ based on $p$ such that $\dim(I, A, M) = \aleph_{\alpha}$.  (This can be proved by an argument similar to the one for the last claim: it suffices to show that if there were an indiscernible sequence $J$ in $M$ extending $I$ such that $\vert J \vert > \aleph_{\alpha}$, then the orthogonality of $p_0$ to all the other $p_i$'s would imply that there is a sequence $J' \subseteq J$ such that $\vert J' \vert < \aleph_{\alpha}$ and $J - J'$ is Morley over $B \cup \bigcup_{i < \omega} I_i$.  This contradicts Fact~\ref{prime_dimensions}.)  Hence (by Lemma~\ref{basechange2}) $\dim(I, M) = \aleph_{\alpha}$.  By the claim above, there is no sequence in $N$ that could be the image of $I$ under an isomorphism, so $M$ and $N$ are not isomorphic.

\end{proof}

\begin{remark}
In the proof of Theorem~\ref{nomadSB}, it is possible for the model $N$ that is constructed to have types of dimension $\kappa$.  For example, consider the case of the theory of an equivalence relation $E$ with infinitely many infinite classes, plus a unary function $S$ in the sort $M / E$ which is a cycle-free bijection.
\end{remark}

\begin{remark}
In conversations with Bradd Hart, he pointed out that it is possible to use a similar construction as in the proof of Theorem~\ref{nomadSB} to prove a slightly more general result: if $T$ is a stable theory with any nomadic type, regular or not, then $T$ does not have the SB~property.  However, in light of Theorem~\ref{nonsstable_SB} and Proposition~\ref{regnomads}, this extra generality is not really necessary.
\end{remark}

\begin{cor}
\label{superstableSB}
If $T$ is superstable and $T$ has a nomadic type then $T$ does not have the SB~property.
\end{cor}

\begin{proof}
By Proposition~\ref{regnomads} and Theorem~\ref{nomadSB}.
\end{proof}

\begin{cor}
\label{class_SB}
If $T$ is a countable theory and $T$ has the SB~property then $T$ is classifiable.
\end{cor}

\begin{proof}
Note that any superstable theory with DOP is multidimensional and thus has a regular nomad by Lemma~\ref{multidimnomad}.  So if we combine Theorem~\ref{OP_NSB}, Theorem~\ref{nonsstable_SB}, and Theorem~\ref{nomadSB}, we see that any theory with the SB~property must be superstable and have NOTOP and NDOP.
\end{proof}

\begin{thm}
\label{nomad_strong_NSB}
If $T$ is a stable theory with a regular nomadic type, then there exists a sequence $\langle M_i : i < \omega \rangle$ of models of $T$ which are pairwise elmentarily bi-embeddable and pairwise nonisomorphic.
\end{thm}

\begin{proof}
We use the same kind of construction as in the proof of Theorem~\ref{nomadSB}.  Let $\alpha, A, B, f$, and $I_i$ be as in the first paragraph of that proof, and for each $i < \omega$, let $M_i$ be $\F^a_\kappa$-prime over $B \cup \bigcup_{j \geq i} I_j$.  The same argument as before shows that these $M_i$'s work.
\end{proof}

Our final result in this section is that nomadic types control when the $\vert T \vert^+$-saturated models of a superstable theory $T$ have the SB~property.  First, we recall a theorem from \cite{bible}:

\begin{fact}
\label{nmd-sat-models}
1. Suppose $T$ is superstable and nonmultidimentional and $M$ is an $\F^a_{\aleph_0}$-saturated model.  Then $M$ is $\F^a_{\aleph_0}$-prime over a set $\bigcup_{s \in S} I_s$, where the $I_s$'s are Morley sequences based on pairwise-orthogonal regular types.

2. Suppose we have $T$, $M$, and $\langle I_s : s \in S \rangle$ as in 1, and $\langle J_s : s \in S \rangle$ is a set of maximal Morley sequences in $M$ such that $\Av(J_s / J_s)$ is nonorthogonal to $\Av(I_s / I_s)$.  Then $M$ is also $\F^a_{\aleph_0}$-prime over $\bigcup_{s \in S} J_s$, and for each $s \in S$, $\vert J_s \vert = \vert I_s \vert$. 
\end{fact}

\begin{proof}
1. This is part of ``Stage~E'' in the proof of Theorem~IX.2.3 of \cite{bible}.

2. By Fact~\ref{regular_dimensions}, $\vert J_s \vert = \vert I_s \vert$.  If $N \prec M$ is an $\F^a_{\aleph_0}$-prime model over $\bigcup_{s \in S} J_s$ and $N$ does not contain some element $a \in I_s$, then the $\F^a_{\aleph_0}$-prime model over $N \cup \left\{a\right\}$ would contain a realization of $\Av(J_s / N)$, contradicting the maximality of $J_s$.

\end{proof}

\begin{thm}
\label{sat-SB-thm}
If $T$ is superstable, then the following are equivalent:

1. $T$ has no nomadic types.

2. The class of $\vert T \vert^{+}$-saturated models of $T$ has the SB~property.

3. For any $\kappa \geq \vert T \vert^{+}$, the class of $\kappa$-saturated models of $T$ has the SB~property.
\end{thm}

\begin{proof}
3 $\Rightarrow$ 2 is trivial, and 2 $\Rightarrow$ 1 follows from Theorem~\ref{nomadSB}.  For 1 $\Rightarrow$ 3, we note that if $T$ has no nomads then $T$ is nonmultidimensional (by Lemma~\ref{multidimnomad}), so by Fact~\ref{nmd-sat-models}, any $\F^a_{\kappa}$-saturated model is  $\F^a_{\kappa}$-prime over a set $\bigcup_{s \in S} I_s$ where the $I_s$'s are Morley sequences based on pairwise-orthogonal regular types.  Now suppose that $M$ and $N$ are bi-embeddable $\F^a_{\kappa}$-saturated models of $T$, witnessed by $f: M \rightarrow N$ and $g: N \rightarrow M$, and $p \in S(\acl^{eq}(\emptyset))$ is a regular type.  Then since $T$ has no nomadic types, there are only finitely many orthogonal conjugates of $p$, say $p_1, \ldots, p_n$.  If $p_i$ is a conjugate of $p$ such that $\dim(p_i, M)$ is minimal and $p_j$ is a conjugate such that $\dim(p_j, N)$ is minimal, then the bi-embeddability of $M$ and $N$ implies that $\dim(p_i, M) = \dim(p_j, N)$.  Continuing by an induction argument, we can show that for any $i$ between $1$ and $n$, $\dim(p_i, M) = \dim(f(p_i), N)$.  This means that if $M$ is $\F^a_{\aleph_0}$-prime over a set $\bigcup_{s \in S} I_s$, where the $I_s$'s are Morley sequences based on pairwise-orthogonal regular types (as in Fact~\ref{nmd-sat-models}), then there is an elementary bijection $h : \bigcup_{s \in S} I_s \rightarrow \bigcup_{s \in S} J_s$, where $J_s$ is a maximal Morley sequence in $N$ based on the type $\Av(f(I) / f(I))$.  By primality and Fact~\ref{nmd-sat-models} again, $h$ can be extended to an isomorphism between $M$ and $N$.

\end{proof}

\subsection{The trivial case}

The word ``trivial'' in the title refers to the standard technical notion of triviality defined below.  The result we prove is that for superstable perfectly trivial theories the Schr\"{o}der-Bernstein property is equivalent to the existence of a nomadic type.  We feel that this provides a little more evidence for Conjecture~\ref{main_SB_conj}, as we will explain more below.

\begin{definition}
(From \cite{trivial_considerations}) 

1.The stable theory $T$ is \emph{trivial} if for any set $A$ and any set $\left\{a, b, c\right\}$ of pairwise $A$-independent elements, $a \ind_{A} b c$.

2. The stable theory $T$ is \emph{perfectly trivial} if for any set $A$ and any $a, b, c$ such that $a \ind_A b$, $a \ind_{Ac} b$.
\end{definition}

\begin{remark}
In the last definition, it does not matter whether we require the elements $a, b,$ and $c$ to be in the home sort of $\mathfrak{C}^{eq}$ or allow them to be imaginary.  If the triviality condition holds for all non-imaginary elements and $a, b, c \in \mathfrak{C}^{eq}$ are pairwise $A$-independent, then we can pick tuples $a', b',$ and $c'$ in the home sort such that $a \in \dcl(a')$, $b \in \dcl(b')$, and $c \in \dcl(c')$, and without loss of generality $a' \ind_{Aa} b c$, $b' \ind_{Ab} a' c$, and $c' \ind_{Ac} a' b'$.  This implies that $a'$, $b'$, and $c"$ are pairwise $A$-independent, so $a' \ind_A b' c'$, and therefore $a' \ind_A b' c'$.
\end{remark}

\begin{lem}
\label{a-orth-triv}
If $T$ is trivial and $p, q \in S(A)$ and $p \bot^a q$ then $p \bot q$.
\end{lem}

\begin{proof}
This follows directly from the definitions: if $a$ realizes $p$, $a \ind_A B$, $b$ realizes $q$, and $b \ind_A B$, then since $a \ind_A b$ (by almost orthogonality), $a \ind_A b B$ (by triviality), so $a \ind_B b$ by monotonicity.
\end{proof}

\begin{thm}
\label{trivial-class}
If $T$ is superstable, perfectly trivial, and nonmultidimensional, $M \models T$, and $J$ is a maximal independent set of elements of $M$ realizing regular types (that is, if $a \in M$ and $\stp(a)$ is regular then $a \nind J$), then $M =\acl(J)$.

\end{thm}

\begin{proof}
This is a direct corollary of Definition~5.1 and Proposition~5.2 of \cite{HPS2}, noticing that nonmultidimensional theories are automatically NDOP.
\end{proof}

\begin{quest}
But can we weaken the perfect triviality hypothesis in the last result to mere triviality?  Note that there are examples of superstable theories that are trivial but not perfectly trivial (see the example discussed at the end of Section~2 of \cite{trivial_considerations}).
\end{quest}

\begin{thm}
\label{trivial-sb-thm}
If $T$ is superstable and perfectly trivial, then $T$ has the SB~property if and only if $T$ has no nomadic types.
\end{thm}

\begin{proof}
Left to right follows directly from Theorem~\ref{nomadSB}.  For right to left, suppose that $M$ and $N$ are bi-embeddable models of $T$, witnessed by $f : M \rightarrow N$ and $g: N \rightarrow M$.  Pick some maximal collection $\langle p_\alpha : \alpha < \mu \rangle$ of pairwise-orthogonal regular types in $S(\acl^{eq}(\emptyset))$ with the following special property: if $h \in \Aut(\mathfrak{C})$ and $h(p_\alpha)$ is nonorthogonal to $p_\beta$, then $p_\alpha = p_\beta$.  From here on, the same argument as in Theorem~\ref{sat-SB-thm} shows that for each $\alpha < \mu$, $\dim(p_\alpha, M) = \dim(f(p_\alpha), N)$ (using Lemma~\ref{a-orth-triv} to show that dimensions match up).  By Theorem~\ref{trivial-class}, $M \cong N$.

\end{proof}

\begin{cor}
If $T$ is trivial and has finite $U$-rank, then $T$ has the SB~property if and only if $T$ has no nomadic types.
\end{cor}

\begin{proof}
By Theorem~\ref{trivial-sb-thm} and Proposition~9 of \cite{trivial_considerations}.
\end{proof}

\bibliography{SB}

\end{document}